%% file: robust_bond_portfolio.tex
\documentclass[12pt]{article}
\usepackage{fullpage,graphicx,psfrag,amsmath,amsfonts,verbatim}
\usepackage[small,bf]{caption}\usepackage{amssymb}
\usepackage{fullpage}
\usepackage{enumitem}
\usepackage{tikz}
\usepackage{listings}
\usepackage{authblk}
\usepackage{nicefrac} % For comparison

\graphicspath{{figures/}}

\DeclareFixedFont{\ttb}{T1}{txtt}{bx}{n}{12} % for bold
\DeclareFixedFont{\ttm}{T1}{txtt}{m}{n}{12}  % for normal

\usepackage{color}
\definecolor{deepblue}{rgb}{0,0,0.5}
\definecolor{deepgreen}{rgb}{0,0.5,0}

\usepackage{hyperref}
\hypersetup{
colorlinks   = true,    % Colours links instead of ugly boxes
urlcolor     = blue,    % Colour for external hyperlinks
linkcolor    = blue,    % Colour of internal links
citecolor    = red      % Colour of citations
}

\input defs.tex
\title{Robust Bond Portfolio Construction\\
via Convex-Concave Saddle Point Optimization}

\author[1]{Eric Luxenberg\footnote{Equal contribution.}}
\newcommand\CoAuthorMark{\footnotemark[\arabic{footnote}]} % get the current value

\author[2]{Philipp Schiele\protect\CoAuthorMark}
\author[1]{Stephen Boyd}
\affil[1]{Department of Electrical Engineering, Stanford University}
\affil[2]{Department of Statistics, Ludwig-Maximilians-Universität München}

\begin{document}
\maketitle

\begin{abstract}
The minimum (worst case) value of a long-only portfolio of bonds, 
over a convex set of yield curves and spreads, can be estimated
by its sensitivities to the points on the yield curve.
We show that sensitivity based estimates are conservative, \ie,
underestimate the worst case value, and that the exact worst case
value can be found by solving a tractable convex optimization problem.
We then show how to construct a long-only bond portfolio that 
includes the worst case value in its objective or as a constraint,
using convex-concave saddle point optimization.
\end{abstract}

\newpage
\tableofcontents
\newpage

\section{Introduction}

We consider a long-only portfolio of bonds, and address
the problem of robust analysis and portfolio construction,
under a worst case framework.  
In this framework we have a set of possible yield curves and bond spreads,
and consider the worst change in value of the portfolio over this 
uncertainty set.

In the analysis problem, considered in \S\ref{s-wc-analysis}, we fix the 
portfolio, and ask what is the worst case change in portfolio value.
We observe that this is a convex optimization problem, readily
solved using standard frameworks or domain specific languages (DSLs)
for convex optimization.
We also consider the linearized version of the same problem,
where the true portfolio value is replaced with its first order 
Taylor approximation.
This approximation can be interpreted as using standard methods to
analyze bond portfolio value using durations.
We show that this is also a convex optimization problem, and is always conservative,
\ie, predicts more of a decrease in portfolio value than the exact method.

In the robust portfolio construction problem,
considered in \S\ref{s-construction}, we seek a portfolio
of bonds that minimizes an objective that includes a robustness term,
\ie, the worst case change in value of the portfolio over
the set of possible yield curves and spreads.
We show that this problem, and its linearized version,
can be formulated as convex-concave saddle point problems,
where we identify the worst case yield and spread and at the same time,
the optimal portfolio.
One interesting ramification of the convex-concave saddle point formulation
is that, unlike in general worst case (minimax) optimization problems,
where there are generally multiple worst case parameters, we need to consider
only one worst case yield curve and set of bond spreads.

In \S\ref{s-dual-method} we show how the convex-concave saddle point
problem can be solved by solving one convex optimization problem.
This is done using the well known technique of expressing the worst case
portfolio value as the optimal value of the dual problem, which converts a 
min-max problem into a min-min problem which we directly solve.
We illustrate the method for a specific case, and in a companion
paper \cite{DSP} explain how the reformulation technique
can be automated, using methods due to 
Juditsky and Nemirovski~\cite{juditsky2021well}.
Using our disciplined saddle point programming framework, we can pose
the robust bond portfolio construction problem in just a few lines 
of simple and natural code, and solve it efficiently.

We present several variations and extensions in \S\ref{s-extensions}, 
including cases where the bond portfolio contains bonds with 
different base yield curves, per-period compounding is used to value
bonds, and a formulation with a robustness constraint as opposed
to an objective term.

Identifying the convex-concave structure of the robust bond portfolio construction
problem is a novel theoretical contribution, and
allows us to use the powerful theory of convex duality to extract insights such
as the existence of a single worst case yield curve, or construct
a robust bond portfolio.  In addition, while this paper was
under review the importance of properly managing risk in a bond portfolio
became quite apparent with the collapse of Silicon Valley Bank, where a
major factor was the bank's exposure to interest rate risk. 

\subsection{Previous and related work}
\paragraph{Bond portfolio construction and analysis.}
Portfolio construction and analysis are well studied problems in finance,
however, most of the literature focuses on equity portfolios. These
approaches often can not be directly applied to bond portfolios, as there are
important differences between the asset classes, such as the
finite maturity of bonds. Yet, by making assumptions about the
re-investment rate, bond portfolios can be constructed and analyzed via
modern portfolio theory (MPT)~\cite{Markowitz1952}. For example, assuming that
the re-investment rate is given by the current spot rate, standard MPT can be
applied to bond portfolios, where the mean and covariance of the bonds can
be derived from sample moments or factor models~\cite{puhle2008bond}.
Similarly,~\cite{korn2006bond} use a factor model for the term structure in an
MPT setting.

Other approaches for bond portfolio construction that are not based on MPT
include exact matching and immunization~\cite{elton2009modern}.
Exact matching is a method for constructing a bond portfolio that minimizes the
required investment amount while ensuring that cash flows arising from
liabilities are being met.
Immunization refers to matching the duration of assets and liabilities,
so that the portfolio value is insensitive to (small) changes in interest
rates. Both of these problems can be formulated as linear programs and so
tractably solved.

Likewise, the factors influencing bond (portfolio) values are well understood,
given the practical implications of the
problem~\cite{fabozzi2012handbook,elton2009modern}.
However, most of the existing literature focuses on parallel shifts in yield
curves and spreads, leading to a trivial worst case scenario. Thus the
literature is sparse when it comes to robust bond portfolio construction as it
relates to possible changes in yield curves and spreads. Instead, most existing
work focuses on the problem of robust portfolio construction under
parameter uncertainty in an MPT framework
(see, \eg,~\cite{tutuncu2004robust,kim2014recent}).

\paragraph{Convex-concave saddle point optimization.}
Convex-concave saddle point problems are a class of optimization problems with
objective functions which are convex in a subset of the optimization variables, and
concave in the remaining variables.  The goal in such problems is to find
a saddle point, \ie, values of the convex variables that minimize
the objective, and values of the concave variables that maximize it.
Convex-concave saddle point problems have
been studied for decades. Indeed, much of the theory of game theory
is based on solving convex-concave saddle point problems, with
early descriptions dating back to the 1920s~\cite{vNeumann1928}, and solutions
based on solving them as a single convex problem via duality dating back to
the 1950s~\cite{morgenstern1953theory}. In their 1983 book, Nemirovski and
Yudin~\cite{nemirovski1983problem} describe the oracle complexity of first
order optimization methods for convex-concave saddle point problems, based on
their previous work on the convergence of the gradient method for
convex-concave saddle point problems~\cite{nemirovski1978cesari}.
Since then, existing work either requires a specific structure of the problem
such as convex and concave variables only being coupled via a bilinear
term~\cite{bredies2016accelerated}, or only under strong assumptions on the
functions' properties~\cite{nemirovski2004prox}.

More recently, Juditsky and Nemirovski~\cite{juditsky2021well}
proposed a general framework for
solving convex-concave saddle point problems with a particular conic structure
as a single convex minimization via dualization.
Building on Juditsky and Nemirovski's work,
the authors of this paper developed
disciplined saddle programming (DSP) \cite{DSP}.
The associated DSL makes it easy to
express a wide class of convex-concave saddle point problems in a natural way;
the problem is then transformed to a single convex optimization problem
using Juditsky and Nemirovski's methods.
This is analogous to disciplined convex programming (DCP),
which makes it easy to specify and solve a wide variety of convex
optimization problems~\cite{diamond-2016}.
The authors' DSP software package allows our formulation of the
robust portfolio construction problem to be specified in just a few lines
of clear code.

\section{Bond portfolio value}
\subsection{Yield curve and spreads}
A bond is a financial contract that obligates the issuer to make a series of 
specified payments over time to the bond holder.
We let $t=1,\ldots, T$ denote time periods, with $t=0$ representing now.
(The periods are usually six months, a typical time between bond coupon payments.)
We represent the bond payments as a vector $c\in \reals_+^T$, where $T$ is the number of
periods, and $\reals_+$ denotes the set of nonnegative  reals.
For each $t = 1, \ldots, T$, $c_t$ is the payment in period $t$ to the bond holder.
A bond has a maturity, which is the period of its last payment; for $t$ larger 
than the maturity, we have $c_t=0$.
The cash flows $c_t$ include coupon payments as well as the payment of 
the face value at maturity.

We consider a portfolio of $n$ bonds, with quantities (also called holdings)
$h=(h_1, \ldots, h_n)\in \reals^n_+$ of each bond,
assumed nonnegative (\ie, only long positions).
We assume that all bonds in the portfolio mature at or before time period $T$.
We first review some basic facts about bonds, for completeness and also
to fix our notation.

Each bond has a known cash flow or sequence of payments,
given by $c_i \in \reals^T_+$, $i=1, \ldots, n$.
We write $c_{i,t}$ to denote the cash flow from bond $i$ in period $t$.
We have $c_{i,t}=0$ for $t$ larger than the maturity of bond $i$.
We let $p \in \reals_+^n$ denote the price of the bonds.
The portfolio value is $V=p^Th$.
The bond prices are modeled using
a base yield curve and spreads for each bond, explained below.

\paragraph{Yield curve.}
The yield curve is denoted by $y \in \reals^T$.
The yield curve gives the discount of a future payment, 
\ie, the current value of a payment of \$1 received in period $t$, denoted $P_t$.
These are given by 
\[
P_t = \exp (-t y_t), \quad t=1, \ldots, T.
\]
We will work with per-period yields, to simplify the formulas,
but following convention, we present all final numerical results as annualized.
(For example, if the periods represent six months, the associated
annualized yields are given by $2y_t$.)
We use continuous compounding for simplicity of notation, but all our results
readily extend to period-wise compounding, where $P_t = (1+y_t)^{-t}$
(see \S\ref{s-extensions}). 

The yield curve gives the
discount of future payments, and captures market expectations with 
respect to macroeconomic factors, fiscal and monetary policy interactions, 
and the vulnerability of private consumption to future (unexpected) shocks.

\paragraph{Bond spreads.}
Bond $i$ has spread $s_i \geq 0$, which means that the bond is priced
at its net present value using the yield curve $y+s_i\ones$, where
$\ones$ is the vector with all entries one.
This is referred to as a `parallel shift' applied to the base yield curve.
We will work with per-period spreads, to simplify the formulas,
but will give final results as annualized.

The spread captures the uncertainty in the cash flow associated with
the bond, such as default or other optionality, which means that 
we value a payment of \$1
from the bond at period $t$ as $\exp(-t(y_t+s_i))$, which is less than or 
equal to $P_t = \exp(-ty_t)$.
The riskier the bond, the larger the spread, which means
future payments are discounted more heavily.

\paragraph{Bond price.}
The price of a bond is modeled as the net present value of its cash flow using
these discounts,
\BEQ \label{eq-bond-price}
p_i = \sum_{t=1}^T c_{i,t} \exp(-t(y_t+s_i)), \quad i=1, \ldots, n.
\EEQ

\paragraph{Portfolio value.}
The portfolio value can be expressed as
\BEQ\label{e-V}
V  =  p^Th = \sum_{i=1}^n\sum_{t=1}^T h_i c_{i,t}\exp(-t(y_t+s_i)).
\EEQ
For reasons mentioned below, it will be convenient to work with the log of 
the portfolio value,
\BEQ\label{e-logV}
\log V  =  \log \left(\sum_{i=1}^n\sum_{t=1}^T h_i c_{i,t}\exp(-t(y_t+s_i))\right).
\EEQ
The portfolio value and log portfolio value are functions of the holdings
$h$, the yield curve $y$, and the spreads $s$, but we suppress this 
dependence to keep the notation light.  (The cash flows $c_{i,t}$ are
fixed and given.)

\paragraph{Convexity properties.}
The portfolio value $V$ is a linear function of $h$,
the vector of holdings, for fixed yield curve and spreads.
If we fix the holdings, $V$ is a convex function of $(y,s)$, the yield
curve and spreads \cite[Chap.~3]{cvxbook}.
The log value $\log V$ is a 
concave function of $h$, for fixed $y$ and $s$,
since it is a concave function of a linear function.
The log value is a convex function
of $(y,s)$, for fixed $h$,
since it can be expressed as
\[
\log V  =  \log \left( \sum_{i=1}^n\sum_{t=1}^T 
\exp(-t(y_t+s_i)+\log h_i +\log c_{i,t}) \right),
\]
which is the log-sum-exp function of an affine function of $(y,s)$
\cite[\S 3.1.5]{cvxbook}.
Thus $\log V$ is a convex-concave function, concave in $h$ and convex in $(y,s)$.

\subsection{Change in bond portfolio value} 
We are interested in the change in portfolio 
value when the yield curve and spreads change
from their current or nominal values $(y^\mathrm{nom},s^\mathrm{nom})$
to the values $(y,s)$, 
with the holdings fixed at $h^\mathrm{nom}$.
We let $V$ denote the portfolio 
value with yield curve and spreads $(y,s)$,
and $V^\mathrm{nom}$ the portfolio value with yield curve and spreads 
$(y^\mathrm{nom},s^\mathrm{nom})$, both with holdings $h^\mathrm{nom}$.
The relative or fractional change in value is given by $V/V^\mathrm{nom}-1$.
It is convenient to work with the change in the log value,
\BEQ\label{e-Delta}
\Delta = \log (V/V^\mathrm{nom}) = \log V - \log V^\mathrm{nom}.
\EEQ
The relative change in value can be expressed in terms of the change
in log value as $\exp \Delta - 1$.
Both $\Delta$ and the relative change in value are readily interpreted.  
For example, $\Delta = -0.15$ means the portfolio
value decreases by the factor $\exp (-0.15) = 0.861$, \ie, a relative
decrease of $13.9\%$.

Since $\log V$ is a convex function of $(y,s)$,
$\Delta$ is also a convex function of $(y,s)$.

\paragraph{First order Taylor approximation.}
The first order Taylor approximation of change in log value, 
denoted $\hat \Delta$, is
\BEQ\label{e-taylor}
\hat \Delta = 
D_\text{yld}^T(y - y^\mathrm{nom}) + D_\text{spr}^T(s - s^\mathrm{nom})
\approx \Delta,
\EEQ
where 
\[
D_\text{yld} = \left.\left(\nabla_{y} \log V\right)
\right\rvert_{y=y^\mathrm{nom},\;s=s^\mathrm{nom}} \in \reals^T, \qquad
D_\text{spr} = \left.\left(\nabla_{s} \log V\right)
\right\rvert_{y=y^\mathrm{nom},\;s=s^\mathrm{nom}} \in \reals^n
\]
are the gradients of the log value with respect to the yield curve and
spreads, respectively, evaluated at the current value 
$(y^\mathrm{nom},s^\mathrm{nom})$.
These are given by
\BEAS
(D_\text{yld})_t &=& -(1/V^\mathrm{nom}) \sum_{i=1}^n
th^\mathrm{nom}_ic_{i,t}\exp(-t(y^\mathrm{nom}_t+s^\mathrm{nom}_i)), \\
(D_\text{spr})_i &=& -(1/V^\mathrm{nom}) \sum_{t=1}^T
th^\mathrm{nom}_ic_{i,t}\exp(-t(y^\mathrm{nom}_t+s^\mathrm{nom}_i)).
\EEAS
(In the first expression we sum over the bonds, while in the second
we sum over the periods.)
The affine approximation \eqref{e-taylor}
is very accurate when $(y,s)$ is near $(y^\mathrm{nom},s^\mathrm{nom})$.

The gradients $D_\text{yld}$ and $D_\text{spr}$ can be given traditional
interpretations.  When $n=1$, \ie, the portfolio consists of a single bond, 
$D_\text{spr}$ is the \emph{duration} of the bond.  When $n=1$ and
$t$ is one of the 12 Treasury spot maturities,
$(D_\text{spr})_t$ is a \emph{key rate duration} of the bond.
(We use the symbol $D$ since the entries of the 
gradients can be interpreted as durations.)
We refer to the Taylor approximation~\eqref{e-taylor}
as the duration based approximation.

\paragraph{A global lower bound.}
Since $\Delta$ is a convex function of $(y,s)$,
its Taylor approximation $\hat \Delta$ 
is a global lower bound on $\Delta$
(see, \eg, \cite[\S 3.1.3]{cvxbook}): For any $(y,s)$ we have
\BEQ\label{e-bnd}
\hat \Delta =
D_\text{yld}^T(y - y^\mathrm{nom}) + D_\text{spr}^T(s - s^\mathrm{nom})
\leq \Delta.
\EEQ
Note that this inequality holds for any $(y,s)$, whereas the
approximation \eqref{e-taylor} is accurate only for $(y,s)$ near
$(y^\mathrm{nom},s^\mathrm{nom})$.
Thus the duration based approximation of the change in log portfolio value
is conservative; the true change in log value will be larger than 
the approximated change in log value.

We can easily obtain a bound on the relative change in portfolio value.
Exponentiating the inequality
\eqref{e-bnd} and using the inequality $\exp u \geq 1+u$, we have
\[
V/V^\mathrm{nom} - 1 = \exp \Delta -1
\geq 
\exp \hat \Delta -1 \geq \hat \Delta.
\]
Therefore $\hat \Delta$ is also a lower bound on the relative change in 
portfolio value using the
duration based approximation; the actual fractional 
change in value will always be more (positive) than the prediction.

\section{Worst case analysis}\label{s-wc-analysis}
In this section we assume the portfolio holdings are known 
and fixed as $h^\mathrm{nom}$,
and consider a nonempty compact convex set $\mathcal U
\subset \reals^T \times \reals^n$ of possible yield curves and spreads.
(We will say more about choices of $\mathcal U$ in \S\ref{s-unc-sets}.)
We define the \emph{worst case portfolio value} as 
\[
V^\mathrm{wc} = \min_{(y,s)\in \mathcal U} V,
\]
\ie, the smallest possible portfolio value over the set of possible
yield curves and spreads.
It will be convenient to work with the worst case 
(\ie, most negative) change in log portfolio value, defined as
\[
\Delta^\mathrm{wc} = \min_{(y,s)\in \mathcal U} 
\Delta= \log V^\mathrm{wc} - \log V^\mathrm{nom}.
\]
When $(y^\mathrm{nom},s^\mathrm{nom}) \in \mathcal U$, 
the worst case log value change is nonpositive.

\subsection{Worst case analysis problem}
We can evaluate $\Delta^\mathrm{wc}$
by solving the convex optimization problem
\BEQ\label{e-wc-prob}
\begin{array}{ll}
\mbox{minimize} & \Delta\\
\mbox{subject to} & (y,s)\in\mathcal{U},
\end{array}
\EEQ
with variables $y$ and $s$.
The optimal value of this problem is $\Delta^\mathrm{wc}$ (from which we 
can obtain $V^\mathrm{wc}$); by solving it, we also find an associated
worst case yield curve and spread, which are themselves interesting.
We refer to \eqref{e-wc-prob} as the \emph{worst case analysis problem}.

\paragraph{Implications.}
One consequence is that we can evaluate $\Delta^\mathrm{wc}$ 
very efficiently using standard methods of convex optimization
\cite{cvxbook}.
Depending on the uncertainty set $\mathcal U$, the problem \eqref{e-wc-prob}
can be expressed very compactly and naturally using
domain specific languages for convex optimization, 
such as CVXPY~\cite{diamond-2016}, CVX~\cite{grant2014cvx},
Convex.jl~\cite{convexjl}, or CVXR~\cite{cvxr2020}.
Appendix~\ref{s-cvxpy-code-wc} gives an example illustrating how simple
and natural the full CVXPY code to solve the worst case analysis problem is.

\paragraph{Maximum element.}
We mention here a special case with a simple analytical solution.
The objective $\Delta$ is 
monotone nonincreasing in its arguments, \ie, increasing any $y_t$ or 
$s_i$ reduces the portfolio value.
It follows that if $\mathcal U$ has a maximum element 
$(y^\mathrm{max},s^\mathrm{max})$, \ie, 
\[
(y^\mathrm{max},s^\mathrm{max}) \geq 
(y,s)~\mbox{for all}~ (y,s)\in \mathcal U
\]
(with the inequality elementwise),
then it is the solution of the worst case analysis problem.
As a simple example, consider
\[
\mathcal U =  \{ (y,u) \mid 
y^\mathrm{min} \leq y \leq y^\mathrm{max},~
s^\mathrm{min} \leq s \leq s^\mathrm{max} \},
\]
\ie, we are given a range of possible values
for each point in the yield curve, and for each spread.
This uncertainty set, which is a hyper-rectangle or box, has maximum element 
$(y^\mathrm{max},s^\mathrm{max})$, which is (obviously) the 
choice that minimizes portfolio value.

More interesting choices of uncertainty sets do not have
a maximum element; for these cases we must numerically solve
the worst case analysis problem \eqref{e-wc-prob}.

\subsection{Linearized worst case analysis problem}
We can replace the objective in \eqref{e-wc-prob} with the lower bound
\eqref{e-bnd} to obtain
the linearized worst case portfolio value problem
\BEQ\label{e-wc-prob-lin}
\begin{array}{ll}
\mbox{minimize} & \hat \Delta=
D_\text{yld}^T(y - y^\mathrm{nom}) + D_\text{spr}^T(s - s^\mathrm{nom}) \\
\mbox{subject to} & (y,s)\in\mathcal{U},
\end{array}  
\EEQ
with variables $y$ and $s$.
Here the objective is affine, whereas in \eqref{e-wc-prob} the 
objective is nonlinear (but convex).
From the inequality \eqref{e-bnd}, solving this linearized
worst case analysis problem gives us a lower bound on $\Delta^\mathrm{wc}$,
as well as a very good approximation when the changes in yield curve and 
spreads are not large.
We refer to \eqref{e-wc-prob-lin} as the \emph{linearized worst case 
analysis problem}, and we denote its optimal value, the 
worst case change in log value predicted by the 
linearized approximation, by $\hat \Delta^\mathrm{wc}$.
This estimate of $\Delta^\mathrm{wc}$ is conservative, \ie,
we have $\hat \Delta^\mathrm{wc} \leq \Delta^\mathrm{wc}$.
The linearized problem is commonly used in practice, and therefore
provides a baseline for comparison.

\subsection{Yield/spread uncertainty sets} \label{s-unc-sets}
In this section we describe some possible choices of 
the yield/spread uncertainty set $\mathcal U$, described as a 
list of constraints.  
Before getting to specifics, we make some comments about high level methods
one might use to construct uncertainty sets.

\paragraph{From all historical data.}  Here we construct $\mathcal U$
from all historical data.  This conservative approach measures the
sensitivity of the portfolio to the yield and spread changing to any
previous value, or to a value consistent with some model of the 
past that we build.

\paragraph{From recent historical data.}  Here we construct $\mathcal U$
from recent historical data, or create a model that places higher weight
on recent data.  The idea here is to model plausible changes to 
the yield curve and spreads using a model based on recent historical values.

\paragraph{From forecasts of future values.}  Here we construct 
$\mathcal U$ as a forecasted set of possible values over the future,
for example a confidence set associated with some predictions.

\paragraph{From current yield and spread estimation error.}
Here $\mathcal U$ represents the set of possible values of the
\emph{current} yield and spreads, which acknowledges that the current 
values are only estimates of some true but unknown value. See, for example,
\cite{filipovic2022stripping} for a discussion of yield curve estimation and a
method that can provide uncertainty quantification.

\subsubsection{Scenarios}
Here $\mathcal U$ is the convex hull of a set of yield curves and spreads,
\[
\mathcal U = \conv \{ (y^1,s^1), \ldots, (y^K,s^K) \},
\]
which is a polyhedron defined by its vertices.
In this case we can think of
$(y^k,s^k)$ as $K$ economic regimes or scenarios.
In the linearized worst case analysis problem, we minimize a linear function
over this polyhedron, so there is always a solution at a vertex, \ie,
the worst case yield curve and spread is one of our scenarios.
In this case we can solve the worst case analysis problem by simply evaluating
the portfolio value for each of our scenarios, and taking the smallest value.
When using the true portfolio value, however, we must solve the problem 
numerically, since the worst case scenario need not be on the vertex;
it can be a convex combination of multiple vertices.

\subsubsection{Confidence ellipsoid}\label{s-ellipsoid}
Another natural uncertainty set is based on a vector Gaussian model
of $(y,s)$, with mean $\mu\in \reals^{T+n}$ and covariance
$\Sigma \in \symm_{++}^{T+n}$, where $\symm_{++}^k$ denotes
the set of symmetric positive definite $k \times k$ matrices.
We take $\mathcal U$ as the associated $(1-\alpha)$-confidence ellipsoid,
\[
\mathcal{U}=\{(y,s)\mid 
((y,s) -\mu)^T \Sigma^{-1}
((y,s) -\mu)\leq F^{-1}(1-\alpha)\},
\]
where $F$ is the cumulative distribution function of a $\chi^2$ distribution
with $T+n$ degrees of freedom.

\subsubsection{Factor model}\label{s-factor-model}
A standard method for describing yield curves and spreads is via a
factor model, with
\[
(y,s) = Zf+v,
\]
where $f \in \reals^k$ is a vector of factors that drive yield curves and
spreads, and $v$ represents idiosyncratic variation, \ie, not due to 
the factors.  (In a statistical model, the entries of $v$ are assumed to 
be uncorrelated to each other and the factor $f$.) 
The matrix $Z \in \reals^{(T+n) \times k}$ gives the 
factor loadings of the yield curve values and spreads.

Typical factors include treasury yields with various maturities,
as well as other economic quantities.
A simple factor model for yields can contain only two or three factors,
which are the first few principal components of historical yield curves,
called level, slope, and curvature~\cite{Litterman1991, cochrane2005bond}.

Using a factor model, we can specify $\mathcal U$ by giving an uncertainty
set $\mathcal F \subset \reals^k$ for the factors, for example as
\[
\mathcal U = \{Zf +v  \mid f \in \mathcal F,~ \|D^{-1} v\|_2^2 \leq 1 \},
\]
where $D$ is a positive diagonal matrix with its entries giving
the idiosyncratic variation of individual yield curve and spread 
values.
We note that while a factor model is typically used to develop 
a statistical model of the yield curve and spreads, we use it here
to define a (deterministic) set of possible values.

\subsubsection{Perturbation description}
The uncertainty set $\mathcal U$ can be described in terms of 
possible perturbations to the current values $(y^\mathrm{nom},s^\mathrm{nom})$.
We describe this for the yield curve only, but similar ideas can be 
used to describe the spreads as well.
We take $y = y^\mathrm{nom}+ \delta$, where $\delta \in \reals^T$ is the 
perturbation to the yield curve.
We might impose constraints on the perturbations such as 
\BEQ\label{e-pert-U}
\delta^\mathrm{min} \leq \delta \leq \delta^\mathrm{max}, \qquad
\sum_{t=1}^T \delta^2_t \leq \kappa, \qquad
\sum_{t=1}^{T-1} (\delta_{t+1}-\delta_t)^2 \leq \omega,
\EEQ
where $\delta^\mathrm{min}, \delta^\mathrm{max}$, $\kappa$, and $\omega$ 
are given parameters, and the first inequalities are elementwise.
The first constraint limits the perturbation in yield for any $t$;
the second limits the mean square perturbation, and the third is a smoothness
constraint, which limits the roughness of the yield curve perturbation.
This is of course just an example; one could add many further constraints, 
such as insisting that the perturbed yield have nonnegative slope,
is concave, or that the perturbation is plausible under a statistical model of
short term changes in yields.

\subsubsection{Constraints}
We can add any convex constraints in our description of $\mathcal U$.
For example, we might add the constraints \eqref{e-pert-U} to a 
factor model, or confidence ellipsoid.
As an example of a constraint related to spreads, we can require
that the spreads are nonincreasing as a function of the bond rating,
\ie, we always have $s_i \leq s_j$ if bond $i$ has a higher rating
than bond $j$.
This is a set of linear inequalities on the vector of spreads $s$.

\section{Robust portfolio construction}\label{s-construction}

In the worst case analysis problem described in \S\ref{s-wc-analysis},
the portfolio is given as $h^\mathrm{nom}$.
Here we consider the case where the portfolio is to be chosen.
We denote the new portfolio as $h$, with 
$h^\mathrm{nom}$ denoting the nominal or current portfolio.
Our goal is to choose $h$, which we do by
minimizing an objective function, subject to some constraints.

\subsection{Nominal portfolio construction problem}
We first describe the nominal bond portfolio construction problem.
We are given a nominal objective function $\phi:\reals^n \to \reals$
which is to be minimized.
The nominal objective function might 
include tracking error against a benchmark,
a risk term, and possibly a transaction cost term if the portfolio is to 
be constructed from the existing portfolio $h^\mathrm{nom}$ 
(see, \eg, \cite{boyd2017multi}).
We will assume that the nominal objective function is convex.

We also have a set of portfolio constraints, which we denote as
$h \in \mathcal H$, where $\mathcal H \subset \reals_+^n$.
The constraint set includes the long-only constraint $h \geq 0$,
as well as a budget constraint, such as $p^T h = p^T h^\mathrm{nom}$,
which states that the new portfolio has
the same value as the original one.  
(This can be extended to take into account transaction costs if needed.)
The constraint set $\mathcal H$ can include constraints on exposures
to regions or sectors, average ratings, duration,
a limit on risk, and so on.
We will assume that $\mathcal H$ is convex.
The nominal portfolio construction problem is
\[
\begin{array}{ll}
\mbox{minimize} & \phi (h)\\
\mbox{subject to} & h\in \mathcal{H},
\end{array}
\]
with variable $h$.  This is a convex optimization problem.

\subsection{Robust portfolio construction problem}\label{s-robust-construction}
To obtain the robust portfolio construction problem we add
one more penalty term to the nominal objective function, which penalizes
the worst case change in value over the given uncertainty set $\mathcal U$.

The term, which we refer to as the \emph{robustness penalty}, is
$-\lambda \Delta^\mathrm{wc}(h)$, where $\lambda>0$ is a parameter
used to trade off the nominal objective $\phi$ and the worst case
change in log portfolio value $-\Delta^\mathrm{wc}(h)$.
Here we write the worst case change in log value with argument $h$,
to show its dependence on $h$.

We arrive at the optimization problem
\BEQ\label{e-wc-synth}
\begin{array}{ll}
\mbox{minimize} & \phi (h) - \lambda \Delta^\mathrm{wc}(h)\\
\mbox{subject to} & h\in \mathcal{H},
\end{array}
\EEQ
with variable $h$.
We refer to this as the \emph{robust bond portfolio construction problem}.
The objective is convex since $\phi$ is convex and $\Delta^\mathrm{wc}$ is
a concave function of $h$.
This means that the robust bond portfolio construction problem is 
convex.

However, the robustness penalty term $-\Delta^\mathrm{wc}(h)$
is not directly amenable to standard convex optimization, 
since it involves a minimization (over $y$ and $s$) itself.  
We will address the question of how to tractably handle the robustness
penalty term below using methods for convex-concave saddle point optimization.

A different framing of the robust bond portfolio construction problem is to
minimize $\phi$ subject to a constraint on $\Delta^\mathrm{wc}$.
This is readily handled, but we defer the discussion to \S\ref{s-extensions}.

\subsection{Linearized robust portfolio construction problem}

As in the worst case analysis problem we can use the 
linearized approximation of the worst case log value instead of the
true log portfolio value, which gives the problem
\BEQ\label{e-wc-synth-lin}
\begin{array}{ll}
\mbox{minimize} & \phi (h) - \lambda \hat \Delta^\mathrm{wc}(h)\\
\mbox{subject to} &h\in \mathcal{H},
\end{array}
\EEQ
with variable $h$.  Here $\hat \Delta^\mathrm{wc}$ is the worst case
change in log portfolio value predicted by the linearized
approximation, \ie, the optimal value of \eqref{e-wc-prob-lin}, as a 
function of $h$.
We note that $\hat \Delta^\mathrm{wc}$ is, like $\Delta^\mathrm{wc}$,
a concave function of $h$.

\subsection{Convex-concave saddle point formulation}
We can write the robust portfolio construction problem \eqref{e-wc-synth} as
\BEQ\label{e-saddle-pt}
\underset{h \in \mathcal H}{\mbox{minimize}} ~
\max_{(y,s) \in \mathcal U} ~
\left( \phi(h) - \lambda \Delta(h, y,s)\right).
\EEQ
(Maximizing $-\lambda \Delta(h,y,s)$ over 
$(y,s)\in \mathcal U$ gives 
$-\lambda \Delta^\mathrm{wc}(h)$.)
The objective in \eqref{e-saddle-pt} is convex in $h$ and 
concave in $(y,u)$, so this is a convex-concave saddle point problem.
Replacing $\Delta$ with $\hat \Delta$ yields the saddle point 
version of the linearized robust portfolio construction problem.

Sion's minimax theorem \cite{sion1958general}
tells us that if $\mathcal{H}$ is compact, when we reverse the order of the 
minimization and maximization we obtain the same value,
which implies that there exists a saddle point $(h^\star,y^\star,s^\star)$, which
satisfies
\[
\phi(h^\star) - \lambda \Delta(h^\star, y,s) \leq 
\phi(h^\star) - \lambda \Delta(h^\star, y^\star, s^\star) \leq 
\phi(h) - \lambda \Delta(h, y^\star, s^\star)
\]
for all $h \in \mathcal H$ and $(y,s)\in \mathcal U$.
The left hand inequality shows that
$\phi(h^\star) - \lambda \Delta(h^\star,y,s)$
is maximized over $(y,s) \in \mathcal U$ by 
$(y^\star,s^\star)$;
the right hand inequality shows that
$\phi(h) - \lambda \Delta(h, y^\star, s^\star)$
is minimized over $h \in \mathcal H$ by $h^\star$.
It follows that
$\phi(h^\star) - \lambda \Delta(h^\star, y^\star, s^\star)$ 
is the optimal value of the robust portfolio construction problem,
$h^\star$ is an optimal portfolio, and $(y^\star,s^\star)$ is 
a worst case yield curve and spread.

\section{Duality based saddle point method}\label{s-dual-method}
The robust bond portfolio construction problem~\eqref{e-saddle-pt} 
is convex, but unfortunately not immediately representable
in a DSL.
In this section we use a well known trick to transform the problem
to one that can be handled directly in a DSL.
Using duality we will express $\Delta^\text{wc}$ as the 
\emph{maximum} of a concave function over
some variables that lie in a convex set.
This method of transforming an inner minimization is not
new; it has been used since the 1950s when Von
Neumann proved the minimax theorem using strong duality in his
work with Morgenstern on game theory~\cite{morgenstern1953theory}.

We now describe the dualization method
for the case when $\mathcal U$ is a polyhedron of the form
\[
\mathcal U = \{ (y,s) \mid A(y,s) \leq b \},
\]
with $A\in \reals^{p \times (T+n)}$ and $b \in \reals^p$,
and the inequality is elementwise, but similar derivations can be 
carried out for other descriptions of $\mathcal U$.

\subsection{Dual form of worst case change in log value}

We assume Slater's condition hold, since any uncertainty set in practice will
have a nonempty relative interior. From strong duality, it follows that
\BEQ\label{e-deltaWC-dual-form}
    \Delta^\text{wc}(h) = \max_{\mu\geq 0,\; \nu}~g(h, \mu,\nu),
\EEQ
where
\BEQ\label{e-Delta-wc-dual}
g(h, \mu,\nu) =\left\{
\begin{array}{ll}
-\log(p^Th)-\mu^Tb-\sum\limits_{i=1}^n\sum\limits_{t=1}^T\zeta(c_{i,t}h_i,\frac{\nu_{i,t}}{t})
&\text{if } A^T\mu -F^T\nu = 0,\ \ones^T\nu = 1,\ \nu \geq 0\\
-\infty & \text{otherwise,}
\end{array}\right.
\EEQ
with
\[
\zeta(x,t) = -t\log(x/t) = t\log(t/x) = t\log(t)-t\log(x),
\]
which is the relative entropy.
Here $\mu \in \reals^p$ and $\nu \in \reals^{nT}$ are variables,
and the matrix
$F\in\reals^{nT\times (T+n)}$ will be defined below.
The relative entropy is convex (see, \eg,~\cite[\S3.2.6]{cvxbook}),
which implies that $g$ is jointly concave in $(h,\mu,\nu)$.
A full derivation of \eqref{e-deltaWC-dual-form} is given in
\S\ref{s-dual-derivation}.

\subsection{Single optimization problem form}\label{s-single-form}
Using \eqref{e-Delta-wc-dual} we can write
the robust bond portfolio construction problem as a single optimization
problem compatible with DSLs,
with variables $h$, $\mu$, and $\nu$:
\[
\begin{array}{ll}
\mbox{minimize} & \phi(h) +\lambda
\left(\mu^Tb+\sum\limits_{i=1}^n\sum\limits_{t=1}^T\zeta(c_{i,t},\frac{\nu_{i,t}}{t})+\log (p^Th)\right)\\
\mbox{subject to} &  \mu \geq 0, \quad A^T\mu - F^T\nu = 0, \quad \nu\geq 0,\quad \ones^T\nu = 1\\
& h\in\mathcal{H}.
\end{array}
\]
This is a convex optimization problem because the objective is convex
and the constraints are linear equality and inequalities.
This form is tractable for DSLs.

\subsection{Automated dualization via conic representation}\label{s-automated}

While for many uncertainty sets explicit dual forms can be derived by hand,
this process can be tedious and error-prone.
In recent work, Juditsky and Nemirovski~\cite{juditsky2021well} present a method for
transforming general structured convex-concave saddle point problems to a single
minimization problem via a generalized conic representation of convex-concave
functions. Similar to  disciplined convex programming (DCP)~\cite{disciplined}, the
method introduces some basic atoms of known convex-concave saddle functions,
as well as a set of rules for combining them
to form composite problems, making it extremely general.

Juditsky and Nemirovski define a general notion of conic representability for
convex-concave saddle point problems
\begin{equation}\label{e-saddle}
\min_{x\in\mathcal{X}}\max_{y\in\mathcal{Y}}\psi(x,y).
\end{equation}
If the convex-concave $\psi$ can be written in this general form, then
\eqref{e-saddle} can be written as a single minimization problem, with variables
comprising $x$ together
with additional variables. See \cite{juditsky2021well} for details on conic
representability, which is beyond the scope of this paper. The set of conically
representable convex-concave functions is large and includes generalized 
inner products of the form $F(x)^TG(y)$ where $F$
is elementwise convex and nonnegative and $G(y)$ is elementwise concave and
nonnegative, and special atoms like weighted log-sum-exp,
$\log(\sum_i x_i\exp(y_i))$, which appears in the robust bond portfolio
construction problem.

%While in \S\ref{s-single-form} we work with a simple uncertainty set (namely,
%polyhedral), this general conic representability approach allows us to perform
%any robust bond portfolio construction problem in which the uncertainty set is
%represented in a DCP compliant fashion and the objective can be expressed
%via a conic representation. Whereas in \S\ref{s-single-form} we had 
%to manually compute the dual of an optimization problem, the advantage of DSP
%is that all of this is fully automated.  

The authors have developed a package for disciplined saddle point programming
called DSP, described in \cite{DSP}.
DSP automates the dualization, following the ideas of Juditsky and Nemirovski,
and allows users to easily express and then solve convex-concave saddle point
problems, including as a special case the
robust bond portfolio construction problem. 
Roughly speaking DSP hides the complexity of dualization from the user,
who expresses the saddle point problem using a natural description.
We refer the reader to \cite{DSP} for a (much) more detailed description of
DSP and its associated domain specific language.

%\clearpage
\subsection{DSP specification}
To illustrate the use of DSP for robust bond portfolio construction,
we give below the code needed to formulate and solve it.
We assume that several objects have already been defined:
\verb|C| is the cash flow matrix,
\verb|H| and \verb|U| are DCP compliant descriptions of 
the portfolio and yield/spread uncertainty set, \verb|phi| is a
DCP compliant convex nominal objective function, and \verb|lamb| is
a positive parameter.

% Minipage to prevent line break
\begin{minipage}{\textwidth}
\begin{lstlisting}[language=Python]
import cvxpy as cp
import dsp

y = cp.Variable(T)
s = cp.Variable(n)
h = cp.Variable(n, nonneg=True)

exponents = []
weights = []
for i in range(n):
    for t in range(T):
        if C[i, t] > 0:
            exponents.append(-(t + 1) * (y[t] + s[i]))
            weights.append(h[i] * C[i, t])

Delta = dsp.weighted_log_sum_exp(cp.hstack(exponents), cp.hstack(weights))

obj = dsp.MinimizeMaximize(phi - lamb * Delta)

constraints = H + U

saddle_problem = dsp.SaddleProblem(obj, constraints)
saddle_problem.solve()

\end{lstlisting}
\end{minipage}

In lines 8--16 we construct an expression for $\Delta$, where 
in line~16 we use
the convex-concave DSP atomic function \verb|weighted_log_sum_exp|.
In line~20 the addition symbol concatenates \verb|H| and \verb|U|, 
which are lists of CVXPY constraints that define $\mathcal H$ and
$\mathcal U$, respectively.
In line~22 we construct the saddle point problem, and in line~23 
we solve it.  The optimal portfolio can then be found in
\verb|h.value|, and the worst case yield and spreads in
\verb|y.value| and \verb|s.value|, respectively.

\section{Variations and extensions}\label{s-extensions}
\paragraph{Periodically compounded growth.}
To handle periodically compounded interest, simply observe that in this case,
\[
y_t = p_t^{-1/t}-1, \quad t=1, \ldots, T,
\]
and the portfolio value is $f(y)=\sum_{t=1}^Tc_t(1+y_t)^{-t}$. This is a convex
function of $y$ because $c_t\geq 0$ and $x^a$ is convex for any $a<0$ and
nonnegative argument.

\paragraph{Multiple reference yield curves.}
We can immediately extend to the case where each bond has its own reference
yield curve $y^i\in\reals^T$ for $i=1,\ldots,n$. This effectively means the
spread can be time varying for each bond. All the convexity properties are
preserved; in fact this just corresponds to an unconstrained $z$ variable in
\S\ref{s-dual-derivation}.

\paragraph{Constrained form.}
It is natural and interpretable to pose~\eqref{e-wc-synth} in constrained form,
that is,
\[
    \begin{array}{ll}
    \mbox{minimize} & \phi (h)\\
    \mbox{subject to} & \Delta^\mathrm{wc}(h)\geq -\eta\\
    & h\in \mathcal{H},
    \end{array}
\]
for some $\eta>0$. For example, one could consider minimizing the tracking error
to a reference bond portfolio, subject to the constraint that the worst case
change in bond portfolio value does not exceed a given tolerance. Using the
conic representability method in \S\ref{s-automated}, we can immediately include
$\Delta^\mathrm{wc}(h)\geq -\eta$ as a DCP compliant constraint. See our
recent paper on DSP~\cite{DSP} for details.

%\clearpage
\section{Examples}

In this section we illustrate worst case analysis and robust portfolio
construction with numerical examples.   The examples all use the data
constructed as described below.
We emphasize that we consider a simplified small problem only so the results
are interpretable, and not due to any limitation in the algorithms used
to carry out worst case analysis or robust portfolio construction, 
which readily scale to much larger problems.

The full source code and data to re-create the results shown here is available
online at \begin{quote}
\url{https://github.com/cvxgrp/robust_bond_portfolio}.
\end{quote}

\subsection{Data}\label{s-data}
We work with a simpler and smaller universe of bonds that is derived from,
and captures the main elements of, a real portfolio.

\paragraph{Bond universe.}
We start with the bonds in the iShares Global Aggregate Bond UCITS ETF 
(AGGG),
which tracks the Bloomberg Global Aggregate Bond
Index~\cite{bloomberg_indices}, a market-cap weighted index
of global investment grade bonds. As of 2022-09-12, AGGG held 10,564 bonds,
which we partition into 20 groups by rating and maturity.
We consider the four ratings AAA, AA, A, and BBB, 
and five buckets of maturities, 0--3, 3--5, 5--10, 10--20,
and 20--30 years.
From each of the 20 rating-maturity groups,
we select the bond in AGGG with
the highest market capitalization.
These 20 bonds constitute the universe we consider.
They are listed in table~\ref{t-bond-universe},
with data as of 2022-09-12.

\begin{table}
    \begin{minipage}{\textwidth} \small
    \centering
\begin{tabular}{lrrrrrr}
Ticker & Rating &  Term to maturity &  Coupon rate &  Distribution & Price \\
 \hline
T 2 \nicefrac{5}{8} 03/31/25     & AAA    & 2.46             & 2.625       & semi-annual  & 101.86 \\
T 1 \nicefrac{1}{4} 12/31/26     & AAA    & 4.21             & 1.250       & semi-annual  & 92.63  \\
T 0 \nicefrac{5}{8} 12/31/27     & AAA    & 5.21             & 0.625       & semi-annual  & 85.43  \\
T 3 \nicefrac{1}{4} 05/15/42     & AAA    & 19.58            & 3.250       & semi-annual  & 128.96 \\
T 3 08/15/52                     & AAA    & 29.84            & 3.000       & semi-annual  & 130.55 \\
FHLMC 0 \nicefrac{3}{8} 09/23/25 & AA     & 2.94             & 0.375       & semi-annual  & 89.13  \\
NSWTC 3 05/20/27                 & AA     & 4.60             & 3.000       & semi-annual  & 106.66 \\
WATC 3 \nicefrac{1}{4} 07/20/28  & AA     & 5.77             & 3.250       & semi-annual  & 110.63 \\
NSWTC 2 \nicefrac{1}{4} 05/07/41 & AA     & 18.56            & 2.250       & semi-annual  & 99.64  \\
BGB 3 \nicefrac{3}{4} 06/22/45   & AA     & 22.69            & 3.750       & annual       & 88.60  \\
JGB 0.4 09/20/25 \#340           & A      & 2.93             & 0.400       & semi-annual  & 88.35  \\
JGB 0.1 09/20/27 \#348           & A      & 4.93             & 0.100       & semi-annual  & 79.62  \\
JGB 0.1 06/20/31 \#363           & A      & 8.68             & 0.100       & semi-annual  & 67.39  \\
JGB 1 12/20/35 \#155             & A      & 13.18            & 1.000       & semi-annual  & 72.75  \\
JGB 1.7 09/20/44 \#44            & A      & 21.94            & 1.700       & semi-annual  & 80.12  \\
SPGB 3.8 04/30/24                & BBB    & 1.54             & 3.800       & annual       & 96.39  \\
SPGB 2.15 10/31/25               & BBB    & 3.05             & 2.150       & annual       & 90.01  \\
SPGB 1.45 10/31/27               & BBB    & 5.05             & 1.450       & annual       & 81.71  \\
SPGB 2.35 07/30/33               & BBB    & 10.79            & 2.350       & annual       & 75.07  \\
SPGB 2.9 10/31/46                & BBB    & 24.05            & 2.900       & annual       & 66.14  \\
\hline
\end{tabular}
\caption{The 20 bond universe used for numerical examples. The prices are
computed as of 2022-09-12.}
\label{t-bond-universe}
\end{minipage}
\end{table}

For each bond we construct its cash flow based on the
coupon rate, maturity, and frequency of coupon payments,
assuming the cash flows are paid at the end of each period.
All bonds in the universe distribute either semi-annual or annual coupons,
so we use a period length of six months. The longest term to maturity
in our universe is 30 years, so we take $T=60$. 
We assemble the cash flows into a matrix $C \in \reals^{20 \times 60}$.
We price the bonds according to~\eqref{eq-bond-price}, using US treasury yield
curve data and spreads that depend on the rating, with data as of
2022-09-12, as listed in table~\ref{t-bond-universe}.

\begin{figure}
    \centering
    \includegraphics[width=0.8\textwidth]{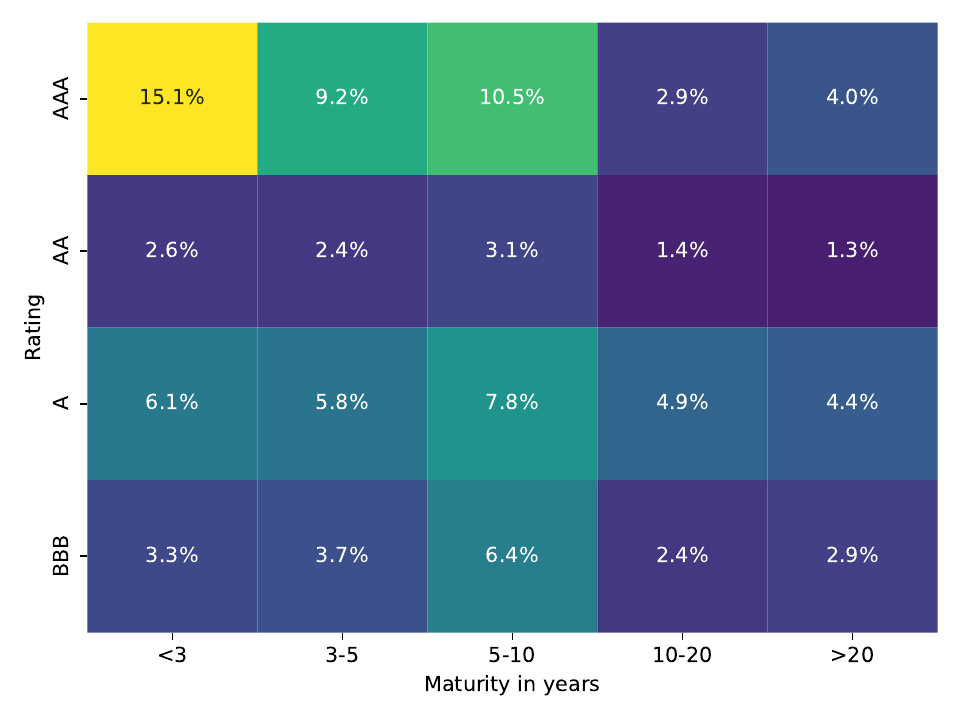}
    \caption{The nominal portfolio weights.}
    \label{f-nominal-portfolio}
\end{figure}

\paragraph{Nominal portfolio.}
Our nominal portfolio $h^\mathrm{nom}$ puts weight on each of our 20 bonds
equal to the total weight of all bonds in the corresponding rating-maturity
group in AGGG.
The weights are shown in figure~\ref{f-nominal-portfolio}.

\subsection{Uncertainty sets}
\begin{figure}
    \centering
    \includegraphics[width=0.9\textwidth]{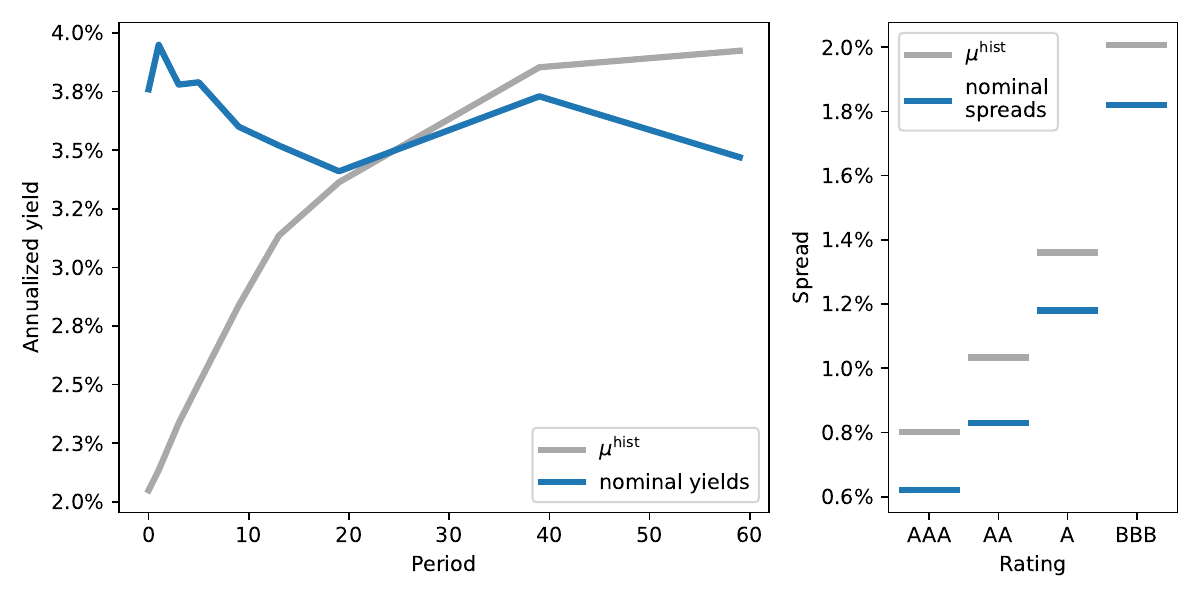}
    \caption{Historical mean values of (annualized) yield and spreads.  The nominal
yield and spreads, date 2022-9-12, are also shown.}
    \label{f-mean}
\end{figure}

We create uncertainty sets using
historical daily yield curves and spreads.
Data for the yield curve is obtained from the US
Treasury~\cite{treasury}, and the spreads are obtained from the
Federal Reserve Bank of St. Louis~\cite{fred},
spanning the period 1997-01-02 to 2022-09-12, for a total of
5,430 observations. As our period length is 6 months, we only consider the 9
key rate durations of 6 months, as well as 1, 2, 3, 5, 7, 10, 20, and
30 years. Joining the 9 rates with the 4 ratings, our total data is represented
as a $5,430 \times 13$ matrix.
The mean value of each column, denoted $\mu^\mathrm{hist}$, as well as the
nominal, \ie, most recent, yields and spreads are shown in figure~\ref{f-mean}.
We use simple linear interpolation to obtain the full 
yield curve $y \in \reals^{60}$.

We model the uncertainty set as a (degenerate) ellipsoid.
We compute the empirical mean and covariance
of the historical data, denoted $\mu^\mathrm{hist}$ and $\Sigma^\mathrm{hist}$,
respectively. We define $Z \in \reals^{80 \times 13}$ as the matrix that maps
the key rates and ratings, $({y^\mathrm{k}}, {s^{r}})\in \reals^{13}$ to the
yields and spreads, $(y,s)\in \reals^{80}$, \ie,
\[
    (y,s) = Z (y^\mathrm{k},s^{r}).
\]
The matrix $Z$ encodes linear interpolation between key rates; other
linear mappings like spline interpolation could also be used.

Our uncertainty set is then defined in terms of key rates and ratings,
\[
\mathcal{U}=\{Z (y^\mathrm{k},s^{r})\mid
\left((y^\mathrm{k}, s^\mathrm{r}) -\mu^\mathrm{hist}\right)^T
\left(\Sigma^\mathrm{hist}\right)^{-1}
\left((y^\mathrm{k},s^\mathrm{r})-\mu^\mathrm{hist}\right)
\leq F^{-1}(1-\alpha)\},
\]
where $F$ is the CDF of a $\chi^2$ random variable with 13 degrees of freedom,
and $\alpha \in (0,1)$ is a confidence level. 
This uncertainty set is a degenerate ellipsoid, with affine dimension $13$.

To represent a modest uncertainty set, we use confidence levels 50\%.
We also consider a more extreme uncertainty set, with confidence level
99\%.
These two uncertainty sets are meant only to illustrate our method; 
in practice, we would likely create uncertainty sets that change over time,
and are based on more recent yield curves and spreads, as opposed to 
a long history of yields and spreads.

\subsection{Worst case analysis}

Table~\ref{t-worst-case} shows worst case change in portfolio 
value for the 50\% and 99\% confidence levels, 
each using both the exact and linearized
methods. The values are given as relative change in portfolio value, \ie,
we have already converted from log returns.
\begin{table}
    \begin{minipage}{\textwidth}
        \centering
        \begin{tabular}{ccc}
            $\alpha$ & Worst case & Linearized \\ \hline
            50\%    & -29.34\%   & -33.43\%   \\
            99\%    & -39.64\%   & -45.95\%   \\
            \hline
        \end{tabular}
        \caption{Worst case change in portfolio value,
for two uncertainty sets, using both the exact and linearized methods.}
        \label{t-worst-case}
    \end{minipage}
\end{table}
For $\alpha = 50\%$,
the exact method gives a worst case change in portfolio value of
around $-29\%$, and the linearized method predicts a change in portfolio value
around four percentage points lower.
For $\alpha=99\%$ the approximation error of the linearized method
is greater, more than six percentage points.
We can see that in both cases the linearized method is conservative,
\ie, predicts a change in value that is lower than the exact method.

\begin{figure}
    \centering
    \includegraphics[width=1.0\textwidth]{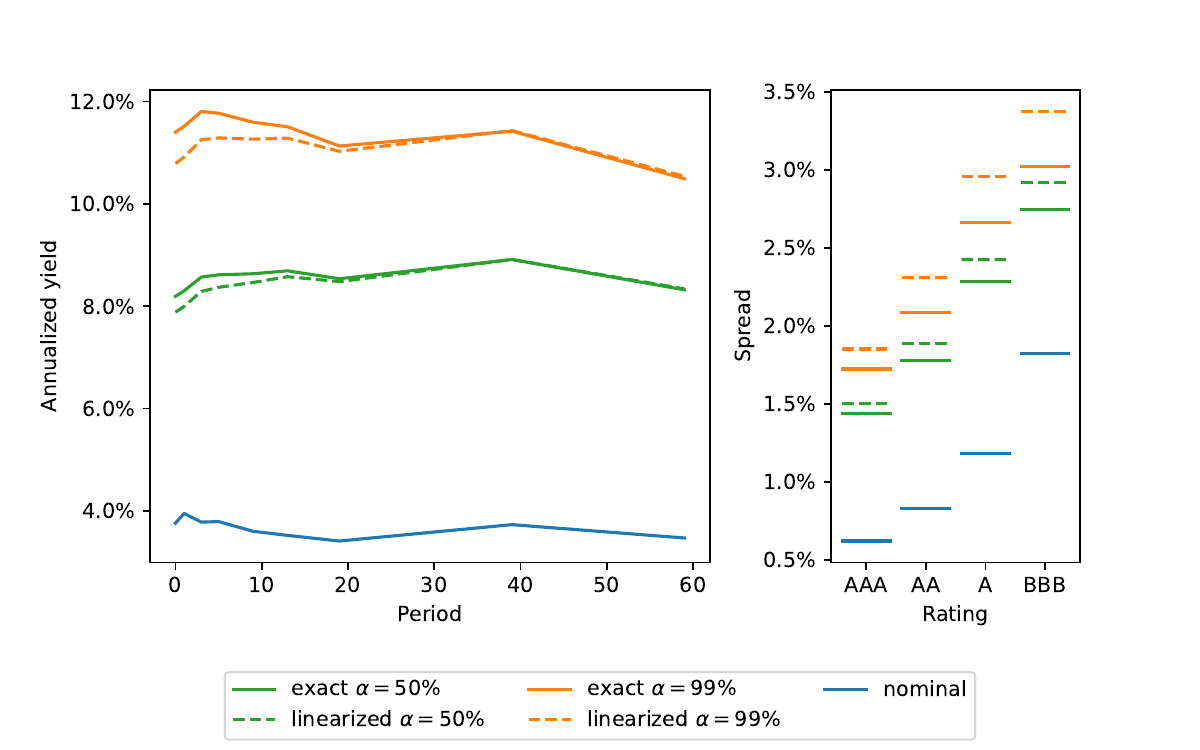}
    \caption{Worst (annualized) yield curve and spreads for the two uncertainty sets,
using the exact and linearized methods.}
    \label{f-worst-case}
\end{figure}

Figure~\ref{f-worst-case} shows the corresponding annualized 
worst case yields and spreads.
For the modest uncertainty case, the linearized and exact methods
produce similar results, while for the extreme case, the linearized
method deviates more for shorter yields and across all spreads.

\subsection{Robust portfolio construction}

We now consider the case where instead of holding the nominal portfolio
$h^\mathrm{nom}$ exactly, we try to track it, with a
penalty term on the worst case change in portfolio value, as described in
\S\ref{s-robust-construction}. Our nominal objective is the turnover distance
between our holdings and the nominal portfolio, given by
\[
\phi(h) = (1/2) \|h-h^{\mathrm{nom}}\|_1.
\]

\begin{figure}
    \centering
    \includegraphics[width=0.8\textwidth]{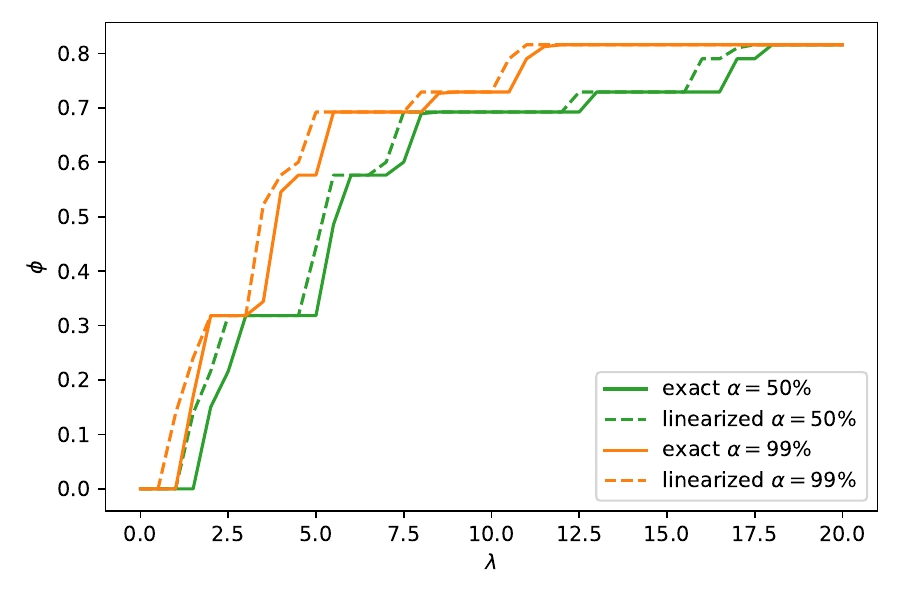}
    \caption{Turnover distance to the nominal portfolio,
        for two uncertainty sets, using both the exact and linearized methods.}
    \label{f-construction}
\end{figure}

Figure~\ref{f-construction} shows the turnover distance for both uncertainty
sets and both the exact and linearized methods. For small values of
$\lambda$, the nominal portfolio is held exactly. As expected, the
turnover distance increases with $\lambda$, but more rapidly so for the
extreme uncertainty set. We also see that the linearized method gives
very similar results. For some values of $\lambda$, the resulting portfolio
obtains the same turnover distance as the exact method, but for other values
it is slightly worse due to the conservative linear approximation.

\begin{figure}
    \centering
    \includegraphics[width=0.79\textwidth]{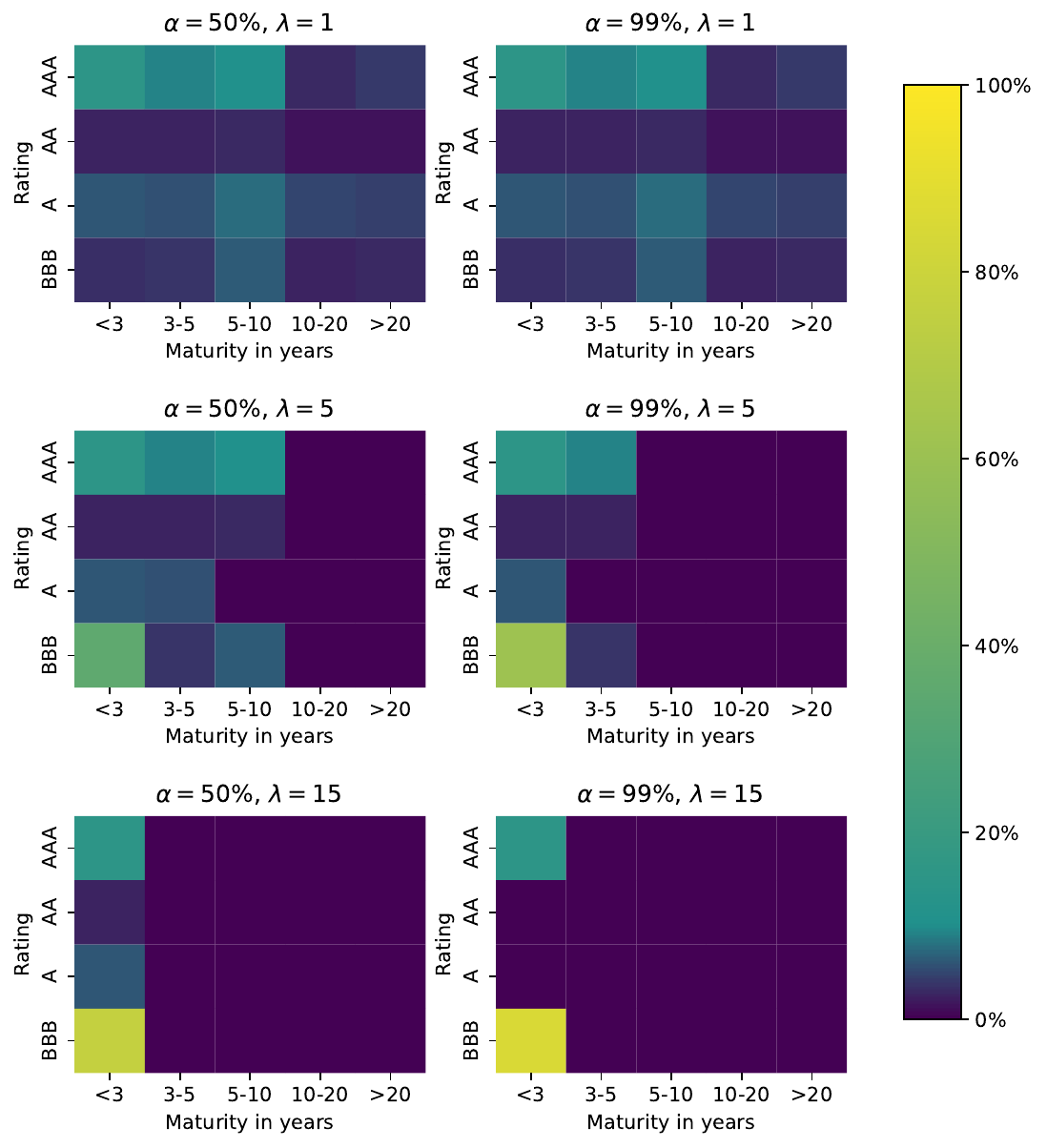}
    \caption{Portfolio holdings, for both uncertainty sets for $\lambda \in \{1,5,15\}$.}
    \label{f-heatmaps}
\end{figure}

Figure~\ref{f-heatmaps} shows the resulting portfolio holdings for the two
uncertainty sets across varying values of $\lambda$. For $\lambda=1$,
the weights are exactly the nominal weights. For $\lambda=5$, the weights have
shifted to shorter maturities, reducing the worst case change in portfolio.
However, this shift is more pronounced for the extreme uncertainty set,
where the weights are only allocated up to the 3--5 year bucket, whereas
for the modest uncertainty set the weights include bonds up to the 5--10 year
bucket. A similar observation can be made for $\lambda=15$, where the
optimization under the modest uncertainty set still allocates in bonds across
all ratings in the bucket containing bonds with less than 3 years to maturity.
In contrast, under the extreme uncertainty set the weights are only allocated
to two bonds in this bucket.
Specifically, the highest weight is assigned to the bond with BBB rating, and
a smaller weight to the AAA bond.
This is explained by the much shorter maturity of the BBB bond
(see table~\ref{t-bond-universe}), which outweighs the lower risk due to
the higher rating of the AAA bond. Indeed, when manually setting the maturities
of these bonds to the same value, we find that weights would be assigned to the
AAA bond for large values of $\lambda$ instead.

\clearpage
\section{Conclusions}
We have observed that the greatest decrease in value of a long-only bond
portfolio, over a given convex set of possible yield curves and spreads,
can be found exactly by solving a tractable convex optimization problem
that can be expressed in just a few lines using a DSL.
Current practice is to estimate the worst case decrease in value using key rate 
durations, which is equivalent to finding the worst case change in
portfolio value using a linearized approximation of the portfolio value.
We show that this estimate is always conservative.  Numerical examples show
that it is a good approximation of the actual worst case value
for modest changes in the yield curve and spread, but less good for 
large uncertainty sets.

We also show that the problem of constructing a long-only bond portfolio
which includes the worst case value over an uncertainty set in its objective
or constraints can be tractably solved by formulating it as a convex-concave
saddle point problem.  Such problems can also be specified in just a few
lines of a DSL.

\section*{Acknowledgements}
We thank Dr. Baruch  Gliksberg for his thoughtful comments and suggestions.

This research was partially supported by ACCESS (AI Chip Center for Emerging 
Smart Systems), sponsored by InnoHK funding, Hong Kong SAR, and by ONR
N000142212121.
P. Schiele is supported by a fellowship within the IFI program of the German
Academic Exchange Service (DAAD).

\clearpage
\bibliography{robust_bond_portfolio.bib}

\clearpage
\appendix

\section{Worst case analysis CVXPY code}\label{s-cvxpy-code-wc}
\begin{lstlisting}[language=Python]
import numpy as np
import cvxpy as cp

y = cp.Variable(T)
s = cp.Variable(n)

V = h @ p

exponents = []
for i in range(n):
    for t_idx in range(T):
        t = t_idx + 1  # account for 0-indexing
        w_it = h[i] * C[i,t_idx]
        if w_it > 0:
            exponents.append(-t * (y[t_idx] + s[i]) + np.log(w_it))

Delta = cp.log_sum_exp(cp.hstack(exponents)) - np.log(V)
obj = cp.Minimize(Delta)
prob = cp.Problem(obj, [A @ cp.hstack([y,s]) <= b])
prob.solve()

\end{lstlisting}

\clearpage
\section{Explicit dual portfolio construction CVXPY code}\label{s-cvxpy-code-dual}
\begin{lstlisting}[language=Python]
import numpy as np
import cvxpy as cp

F_1 = np.tile(np.eye(T), (n, 1))
F_2 = np.repeat(np.eye(n), repeats=T, axis=0)
F = np.hstack([F_1, F_2])

lam = cp.Variable(len(b), nonneg=True)
nu = cp.Variable(n * T, nonneg=True)

h = cp.Variable(n, nonneg=True)

B = 1

term = 0
for i in range(n):
    for t in range(1, T):
        nu_i_t = nu[i * T + t]
        term -= cp.rel_entr(C[i, t] * h[i], nu_i_t / t)

obj = cp.Maximize(-lam @ b + term - np.log(B))
constraints = [
    A.T @ lam == F.T @ nu,
    cp.sum(nu) == 1,
    p @ h == B,
]
prob = cp.Problem(obj, constraints)
prob.solve()

\end{lstlisting}

\clearpage
\section{Derivation of dual form}\label{s-dual-derivation}

We now derive the dual form of the worst case portfolio construction problem
for the case where the uncertainty set in polyhedral.
We note that the worst case $\log$ change in portfolio value for a fixed
$h$,
$\Delta^\text{wc}(h)$
is given by the optimal value of the optimization problem
\[
\begin{array}{ll}
\mbox{minimize} & \log\left(\sum_{i=1}^n\sum_{t=1}^Th_ic_{i,t}\exp(-t(y_t+s_i))\right) -\log(p^Th)\\
\mbox{subject to} & A(y,s)\preceq b.
\end{array}
\]
We have written this problem with $(y,s)$ explicitly, instead of
with $x$, to emphasize the objective's dependence on each component. We note
that due to our budget constraint, $ \log(V(y,s))=\log(p^Th)$.

In order to obtain a closed form dual, we introduce a new variable
$z\in\reals^{nT}$, where we think of $z_{i,t}$ as corresponding to $
y_t+ s_i$.
This is a very general formulation which allows each
bond to be associated with its own yield curve $y_i\in \reals^T$, with $z_{i,t}$ corresponding to the $t$'th entry of the $i$'th
bond's yield curve.
Since we model each bond as having its own yield curve, this formulation
generalizes the earlier treatment with yields and spreads handled separately. We
can recover the original structure with the linear constraints
\BEQ\label{e-F-mat}
z_{i,t} = y_t+s_i,\qquad t=1,\ldots,T,\quad i=1,\ldots,n,
\EEQ
which are representable as $z=Fx$ for an appropriate
$F\in\reals^{nT\times (T+n)}$.
As such, the problem is equivalent to
\BEQ\label{e-wc-prob-full}
\begin{array}{ll}
\mbox{minimize} & \log\left(\sum_{i=1}^n\sum_{t=1}^Th_ic_{i,t}\exp(-tz_{i,t})\right)-\log(p^Th) \\
\mbox{subject to} & Ax \preceq b,\quad z=Fx.
\end{array}
\EEQ
Strong duality tells us that $\Delta^\text{wc}(h)$ is equal to the optimal value of the
dual problem of \eqref{e-wc-prob-full} \cite[\S 5.2]{cvxbook}.

\paragraph{Dual problem.}
We derive the dual of the problem
\BEQ
\begin{array}{ll}
\mbox{minimize} & \log\left(\sum_{i=1}^n\sum_{t=1}^Th_ic_{i,t}\exp(-tz_{i,t})\right)-\log(p^Th) \\
\mbox{subject to} & Ax \preceq b,\quad z=Fx.
\end{array}
\EEQ
First, with $f$ the log-sum-exp function $f(x)=\log\left(\sum \exp x_i\right)$, we observe that our problem can be rewritten as
\[
\begin{array}{ll}
\mbox{minimize} & f(Cz+d)-\log(p^Th) \\
\mbox{subject to} & Ax \preceq b,\quad z=Fx.
\end{array}
\]
We define $C$ to be the diagonal matrix with $C_{i,t}=-t$, where we are
using unwound vectorized indexing for $z$, and $d\in\reals^{nT}$ to be the
vector with $d_{i,t} = \log(c_{i,t}h_i)$.
Then, the Lagrangian is given by
\[
L(z,x,\mu,\nu) = f(Cz+d) +\mu^T(Ax-c)+\nu^T(z-Fx)-\log(p^Th).
\]
The Lagrange dual function is given by
\[
g(\mu,\nu) =\inf_{z,x} L(z,x,\mu,\nu) = \inf_z\left(f(Cz+d)+\nu^Tz\right)
+\inf_x\left(\mu^TAx-\nu^TFx\right) -\mu^Tc-\log(p^Th).
\]
The second term is equal to $-\infty$ unless $A^T\mu+F^T\nu=0$, so
this condition will implicitly restrict the domain of $g$.
Now, note that with $g(z)=f(Cz+d)$, the first term can be rewritten as
\BEAS
\inf_z\left(f(Cz+d)+\nu^Tz\right) &=& \inf_z\left(g(z)+\nu^Tz\right)\\
&=& -\max_z \left( -\nu^Tz-g(z)\right)\\
&=& -g^*(-\nu),
\EEAS
where $g^*(y)=\max{y}y^Tz-g(z)$ is the conjugate of $g$.
\cite[\S3.3.1]{cvxbook}.

We now use two facts from \cite[\S3.3.2]{cvxbook}. First, in general the conjugate of the linear
precomposition $\phi(z)=\rho(Cz+d)$ can be written in terms of the conjugate of
$\rho$ as $\phi^*(y) = \rho^*(C^{-T}y)-d^TC^{-T}y$.
Second, the dual of the log-sum-exp function $f$ is
\[
f^*(y) = \left\{\begin{array}{ll}
\sum_iy_i\log(y_i)    & \text{ if } y\geq 0,\quad \ones^Ty = 1\\
\infty    & \text{ otherwise.}
\end{array}\right.
\]
Combining these two facts, and expanding terms, we find that
\[
g(\mu,\nu) =\left\{
\begin{array}{ll}
-\log(p^Th)-\mu^Tb-\sum\limits_{i=1}^n\sum\limits_{t=1}^T\zeta(c_{i,t}h_i,\frac{\nu_{i,t}}{t})
&\text{if } A^T\mu -F^T\nu = 0,\ \ones^T\nu = 1,\ \nu \geq 0\\
-\infty & \text{otherwise,}
\end{array}\right.
\]
with
\[
\zeta(x,t) = -t\log(x/t) = t\log(t/x) = t\log(t)-t\log(x).
\]

Thus the robust bond portfolio problem can be written as
\[
\begin{array}{ll}
\mbox{minimize} & \phi(h) -\lambda \max\limits_{\mu,\nu}~g(\mu,\nu)\\
\mbox{subject to} & \mu\geq 0,\quad A^T\mu -F^T\nu = 0,\quad \ones^T\nu = 1,\quad \nu \geq 0\\
& h\in \mathcal{H},\\
\end{array}
\]
where we have moved the implicit constraints in the definition of $g$ to
explicit constraints in the optimization problem. Note this equivalent
optimization problem has new variables $\mu$ and $\nu$. By using that
$-\lambda\max_{\mu,\nu}g(\mu,\nu)= \min_{\mu,\nu}-\lambda g(\mu,\nu)$ and
collecting the minimization over $h$, $\mu$, and $\nu$, we obtain
the form in \S\ref{s-single-form}.

\end{document}

%% file: defs.tex
\newcommand{\ones}{\mathbf 1}
\newcommand{\reals}{{\mbox{\bf R}}}

\newcommand{\symm}{{\mbox{\bf S}}}  % symmetric matrices

\newcommand{\conv}{\mathop{\bf conv}}

\newcommand{\eg}{{\it e.g.}}
\newcommand{\ie}{{\it i.e.}}

\newcommand{\BEAS}{\begin{eqnarray*}}
\newcommand{\EEAS}{\end{eqnarray*}}
\newcommand{\BEA}{\begin{eqnarray}}
\newcommand{\EEA}{\end{eqnarray}}
\newcommand{\BEQ}{\begin{equation}}
\newcommand{\EEQ}{\end{equation}}
\newcommand{\BIT}{\begin{itemize}}
\newcommand{\EIT}{\end{itemize}}

\newcounter{algorithmctr}
\renewcommand{\thealgorithmctr}{\arabic{algorithmctr}}
   {\mbox{}\\*[\parskip]\begin{minipage}{\linewidth}%
       \refstepcounter{algorithmctr}\begin{list}{}{%
       \setlength{\rightmargin}{0\linewidth}%
       \setlength{\leftmargin}{.05\linewidth}}%
       \rmfamily\small
       \item[]{\setlength{\parskip}{0ex}\hrulefill\par%
        \nopagebreak{\bfseries\textsf{Algorithm \thealgorithmctr~}}}}%
   {{\setlength{\parskip}{-1ex}\nopagebreak\par\hrulefill\\*[2ex]\par}% 
   \end{list}\end{minipage}}